\documentclass[12pt]{amsart}
\input{xypic}
\usepackage{amssymb}

\newtheorem{theorem}{Theorem}[section]
\newtheorem{lemma}[theorem]{Lemma}

\newtheorem{example}[theorem]{Example}
\newtheorem{proposition}[theorem]{Proposition}

\newtheorem{remark}[theorem]{Remark}

\numberwithin{equation}{section}

\newcommand{\mono}{\hookrightarrow}

\newcommand{\epi}{\mbox{$\to$\hspace{-0.35cm}$\to$}}

\def\umono{\ar@{_{(}->}[u]}
\def\uumono{\ar@{_{(}->}[uu]}

\def\lmono{\ar@{_{(}->}[l]}
\def\llmono{\ar@{_{(}->}[ll]}

\newcommand{\Hom}{{\rm Hom \,}}
\newcommand{\map}{{\rm Map}\,}

\newcommand{\Z}{{\mathbb Z}}

\title[Cellularization and fusion]{Cellularization of classifying spaces and fusion properties of finite groups}

\author{Ram\'on J. Flores and J\'er\^ome Scherer}
\thanks{Both authors are supported by MEC grant MTM2004-06686
and the second author by the Ram\'on y Cajal program, MEC, Spain.}

\date{\today}

\begin{document}

\begin{abstract}
One way to understand the mod $p$ homotopy theory of classifying
spaces of finite groups is to compute their
$B\Z/p$-cellularization. In the easiest cases this is a
classifying space of a finite group (always a finite $p$-group).
If not, we show that it has infinitely many non-trivial homotopy
groups. Moreover they are either $p$-torsion free or else
infinitely many of them contain $p$-torsion. By means of
techniques related to fusion systems we exhibit concrete examples
where $p$-torsion appears.
\end{abstract}

\maketitle

\section*{Introduction}
Let $A$ be a pointed space. Although the idea of building a space
as a homotopy colimit of copies of another fixed space $A$ goes
back to Adams (\cite{Adams78}) in the framework of the
classification of the acyclics of a certain generalized cohomology
theory, it was in the early 90's that Dror-Farjoun
(\cite{Farjoun95}) and Chach\'olski (\cite{MR97i:55023})
formalized and developed this idea in the wider context of
cellularity classes. Thus a space will be called \emph{A-cellular}
if it can be built from $A$ by means of pointed homotopy colimits.
There exists an $A$-cellularization functor $CW_A$ that provides
the best possible $A$-cellular approximation: the natural map
$CW_A X \rightarrow X$ induces a weak equivalence between pointed
mapping spaces $\map_*(A, CW_A X) \simeq \map_*(A, X)$.

Our interest lies essentially in using the cellularization functor
to study the $p$-primary part of the homotopy of the classifying
space of a finite group $G$. This approach was already suggested
by Dror-Farjoun in \cite[Example 3.C.9]{Farjoun95}, and proved to
be very fruitful in the last years; we can remark work of
Bousfield (\cite{Bousfield97}) which describes cellularization of
nilpotent spaces with regard to Moore spaces $M(\Z /p,n)$, or the
relationship recently discovered (\cite{MR2003a:55024}) between
the $B\Z /p$-cellularization of spaces and the $\Z
/p$-cellularization in the category of groups.

In the present paper we focus our attention in describing the
$B\Z/p$-cellularization of $BG$, where $G$ is a finite group. From
the above description of the cellularization functor, we know this
is a space which can be built from $B\Z/p$ by means of push-outs,
wedges, and telescopes, and which encodes all the information
about the pointed mapping space $\map_*(B\Z/p, BG)$, which is
isomorphic to $\Hom(\Z/p, G)$. A first attempt to understand
$CW_{B\Z /p}BG$ was undertaken in \cite{Ramon}, where in
particular the fundamental group of the cellularization was
described as an extension of the group-theoretical $\Z
/p$-cellularization of a subgroup of $G$ by a finite $p$-torsion
free abelian group (see Sction 1 for details). In the present
note we concentrate thus on the higher homotopy
groups $\pi_n(CW_{B\Z /p}BG)$,
for $n\geq 2$.

Our first result establishes that the homotopy of $CW_{B\Z /p}BG$
is subject to a dichotomy finite-infinite. This contrasts with the
results in ~\cite{CCS2}: for easier spaces to work with, such as
$H$-spaces or classifying spaces of nilpotent groups, the
assumption that the mapping space $\map_*(B\Z/p, X)$ is discrete
implies that $CW_{B\Z/p} X$ is aspherical.

\vspace{.5cm}

\noindent {\bf Theorem \ref{dichotomy}.} \emph{Let $G$ be a finite
group. Then $CW_{B\Z/p} BG$ is either the classifying space of a
finite $p$-group, or it has infinitely many non-trivial homotopy
groups.}

\vspace{.5cm}

An important role in the proof of this result is played by Levi's
dichotomy theorem \cite[Theorem 1.1.4]{MR96c:55019} about the
homotopy structure of the Bousfield-Kan $p$-completion of $BG$,
and just like this, it is very much in the spirit of the
McGibbon-Neisendorfer theorem (\cite{Neisendorfer84}) and other
similar results that have arisen in the last years in the context
of localization (\cite{Casacuberta93}, \cite{Bastardas02}). All of
these dichotomies share the common feature that they essentially
concern objects whose homotopy, if infinite, is a $p$-torsion
invariant. In the case of $B\Z /p$-cellularization, however, it
cannot be deduced from the construction of the functor that the
homotopy groups of a $B\Z /p$-cellular space should be necessarily
$p$-groups. Hence it becomes evident that one should start a more
precise study of the higher homotopy of $CW_{B\Z /p}BG$. In fact
Chach\'olski's description of the cellularization as a small
variation of the homotopy fiber $\bar P_{B\Z/p}$ of the
$B\Z/p$-nullification map makes us investigate the difference
between these two functors.

\vspace{.5cm}

\noindent {\bf Theorem~\ref{step}.} \emph{Let $G$ be a finite
group. Then either the cellularization $CW_{B\Z/p} BG$ has
infinitely many homotopy groups containing $p$-torsion or it fits
in a fibration
$$
CW_{B\Z/p} BG \longrightarrow BG\longrightarrow \prod_{q\neq p}
(BG)^{\wedge}_q ,
$$
where the (finite) product is taken over all primes $q$ dividing
the order of $G$, and the right map is the product of the
completions.}

\vspace{.5cm}

Note that in the second case the upper homotopy of the
cellularization is that of $\Omega  \prod_{q\neq p}
(BG)^{\wedge}_q$. In particular,
if $G$ is $\Z /p$-cellular, the results  \cite[20.10]{MR97i:55023}
and \cite[3.5]{Ramon} imply that $CW_{B\Z/p} BG$ is different
from $\bar P_{B\Z/p}$ exactly when there appears (a lot of)
$p$-torsion in the higher homotopy groups. Curiously enough this
yields many examples of $B\Z /p$-cellular spaces whose higher
homotopy is $p$-torsion \emph{free}, including the cellularization
of classifying spaces of all symmetric groups. The main question
which remains unanswered up to this point is thus to determine if
there actually exist groups for which $CW_{B\Z/p} BG$ does contain
$p$-torsion in its higher homotopy groups! This turns out to
depend heavily on the $p$-complete classifying space
$BG^{\wedge}_p$.

\vspace{.5cm}

\noindent {\bf Proposition \ref{completion}.} \emph{Let $G$ be a
group generated by order $p$ elements which is not a $p$-group.
Then the universal cover of $CW_{B\Z/p} BG$ is $p$-torsion free if
and only if the $p$-completion of $BG$ is $B\Z/p$-cellular.}

\vspace{.5cm}

We remark (see Proposition~\ref{nontrivialmap}) that specific
representations of $G$ in some $p$-completed compact Lie group
detect the presence of $p$-torsion in the higher homotopy groups
of the cellularization of $BG$. This observation definitively
shifts the problem to the study of fusion systems, a notion which
has recently led Broto, Levi, and Oliver to the concept of
$p$-local finite groups in a topological context (see
\cite{Broto03} and \cite{Broto032}). The map $BG^{\wedge}_p
\rightarrow BK^{\wedge}_p$ we are looking for --where $K$ is some
compact Lie group-- is indeed best understood as a fusion
preserving representation of the Sylow $p$-subgroup of $G$ in $K$
(this relationship is well explained by Jackson
in~\cite{Jackson04}), and turns out
to be trivial when precomposing
with any map $B\Z /p\longrightarrow BG^{\wedge}_p$. To our
knowledge, the kind of representations we are looking for had not
yet been described
in the extensive literature devoted to the classification of
representations of finite groups.

We prove that the Suzuki group $Sz(8)$ admits such a
representation into $U(7)$ at the prime $2$, and so there exist infinitely many
homotopy groups of $CW_{B\Z/2} BSz(8)$ containing $2$-torsion. The
same phenomenon actually occurs for all Suzuki groups, and
moreover the tools used can be applied to other groups at odd
primes (see Example~\ref{oddprimes}). Then, we can establish

\vspace{.5cm}

\noindent {\bf Theorem \ref{SuzukiCW}.} \emph{For every
integer $n$ of the form $2^{2k+1}$ the
$B\Z /2$-cellularization of $BSz(n)$ has 2-torsion in an infinite
number of homotopy groups. Likewise, if $p$ is an odd prime and
$q$ is any integer of the form $mp^k+1$ with $k \geq 2$, then
the $B\Z /p$-cellularization of  $BPSL(q)$ has $p$- torsion in an infinite
number of homotopy groups.}

\vspace{.5cm}

We point out that the representation $BSz(8) \rightarrow
BU(7)^{\wedge}_2$ cannot be induced by a homomorphism of groups
$Sz(8) \rightarrow U(7)$, not even composed with an Adams
operation. Thus our example could be compared with the map
$BM_{12}\to BG_2$ constructed by Benson and Wilkerson in
\cite{MR1349127}. It is also related with work of Mislin-Thomas
\cite{MR985537}, and more recently of Broto-M{\o}ller \cite{BM04}.

\medskip

{\bf Acknowledgements}. We would like to thank Wojciech
Chach\'olski, Jesper Grodal, Ran Levi, Michel Matthey, and Guido
Mislin for helpful discussions during the Arolla conference in
August 2004. On our way to Theorem~\ref{Suzuki} we were helped by
Jacques Th\'evenaz, who constructed for us the group in
Example~\ref{Thevenaz},  and Albert Ruiz and Jesper M{\o}ller, with
whom we computed a few fusion systems. We warmly thank
Carles Broto for many valuable suggestions,
Bob Oliver, who hinted at the Suzuki groups
for their beautiful fusion properties. and Antonio
Viruel for recognizing the same features in the projective special
linear groups.

\section{Background}
There are two main ingredients used to analyze the cellularization
$CW_{B\Z/p} BG$. The first one is a general recipe, the second a
specific simplification in the setting of classifying spaces of
finite groups.

One of the most efficient tools to compute cellularization
functors is Chach\'olski's construction of $CW_A$ out of the
nullification functor $P_A$, see \cite{MR97i:55023}. His main
theorem states that $CW_A X$ can be constructed in two steps.
First consider the evaluation map $\vee_{[A, X]_*} A \rightarrow
X$ (the wedge is taken over representatives of all pointed
homotopy classes of maps) and let $C$ denote its homotopy cofiber.
Then $CW_A X$ is the homotopy fiber of the composite map $X
\rightarrow C \rightarrow P_{\Sigma A} C$. We will actually need a
small variation of Chach\'olski's description of the
cellularization.

\begin{lemma}
\label{betterCW}
Let $A$ and $X$ be pointed connected spaces and choose a pointed
map $A \rightarrow X$ representing each unpointed homotopy class
in $[A, X]$. Denote by $D$ the homotopy cofiber of the evaluation
map $ev: \bigvee_{[A, X]} A \rightarrow X$. Then $CW_A X$ is
weakly equivalent to the homotopy fiber of the composite $X
\rightarrow D \rightarrow P_{\Sigma A} D$.
\end{lemma}

\begin{proof}
According to \cite[Theorem~20.3]{MR97i:55023} one has to check
that composing any map $f: A \rightarrow X$ with $ev: X
\rightarrow D$ is null-homotopic. But such a map $f$ is freely
homotopic to one in $[A, X]$. Thus $ev \circ f$ is freely
homotopic to the constant map, and so it must also be
null-homotopic in the pointed category.
\end{proof}

In \cite{Ramon}, the first author focused on the situation when $A
= B\Z/p$ and $X$ is the classifying space of a finite group $G$. A
first reduction can always be done by considering the subgroup
$\Omega_1(G)_p$ of $G$ generated by the elements of order $p$.
This notation is the standard one in group theory, see for example
\cite{Gor80}. In \cite{Ramon} and \cite{MR2003a:55024} the
terminology ``socle" and the corresponding notation $S_{\Z/p} G$
was used instead. We have indeed an equivalence of pointed mapping
spaces $\textrm{Map}_*(B\Z/p, BG) \simeq \textrm{Map}_*(B\Z/p,
B\Omega_1(G)_p)$, which means that $B\Omega_1(G)_p \rightarrow BG$
is a $B\Z/p$-cellular equivalence. A second simplification
consists then in computing the group theoretical cellularization
$CW_{\Z/p} G \cong CW_{\Z/p} \Omega_1(G)_p$, which is known to be
a (finite) central extension of $\Omega_1(G)_p$ by a group of
order coprime with $p$. Since $B CW_{\Z/p} G \rightarrow BG$ is a
$B\Z/p$-cellular equivalence, we will from now assume that $G$ is
a finite $\Z/p$-cellular group (it can be constructed out of the
cyclic group $\Z/p$ by iterated colimits). Examples of such groups
are provided at the prime $2$ by dihedral groups, Coxeter groups,
and in general by finite $p$-groups generated by order $p$
elements.

The first author described the fundamental group of $CW_{B\Z/p}
BG$:

\begin{proposition}\cite[Theorem 4.14]{Ramon}
\label{fundamentalgroup}
Let $G$ be a finite $\Z/p$-cellular group. Then the fundamental group
$\pi = \pi_1 CW_{B\Z/p} BG$ is described as an extension $ H \mono
\pi \epi G$  of $G$ by a finite $p$-torsion free abelian group
$H$.
\end{proposition}

He also showed that $BG$ is $B\Z/p$-cellular if and only if $G$ is
a finite $\Z/p$-cellular $p$-group. This paper is an attempt to
understand $CW_{B\Z/p} BG$ when $G$ is not a $p$-group.

\section{A cellular dichotomy}

We establish in this section a first dichotomy result. The space
$CW_{B\Z/p} BG$ has infinitely many non-trivial homotopy groups
unless it is aspherical. We will say more about the higher
homotopy groups in the next section.

\begin{theorem}
\label{dichotomy}
Let $G$ be a finite group. Then $CW_{B\Z/p} BG$ is either the
classifying space of a finite $p$-group, or it has infinitely many
non-trivial homotopy groups.
\end{theorem}

\begin{proof}
If $CW_{B\Z/p} BG$ is the classifying space of a group, then it
must be a finite $p$-group as shown by the first author in
\cite[Theorem 4.14]{Ramon}. Assume therefore that it is not so,
i.e. there exists a prime $q$ different from $p$ dividing the
order of $G$. Recall that we may always suppose that $G$ is
$\Z/p$-cellular. Therefore the evaluation map $*\Z/p \rightarrow
G$ taken over all homomorphisms $\Z/p \rightarrow G$ is
surjective. Hence the homotopy cofiber $C$ of the evaluation map
$\vee_{[B\Z/p, BG]} B\Z/p \rightarrow BG$, as defined above, is a
simply-connected space and so is its nullification $P_{\Sigma
B\Z/p} C$.

The fibration from \cite{MR97i:55023}
$$
CW_{B\Z/p} BG \longrightarrow BG \longrightarrow P_{\Sigma B\Z/p}
C
$$
shows that the cellularization shares many homotopical properties
with the nullification $P_{\Sigma B\Z/p} C$. We want to analyze
the $q$-completion of this last space. But since $B\Z/p$ is
$H\Z/q$-acyclic, we see that the composite
$$
BG \longrightarrow C \longrightarrow P_{\Sigma B\Z/p} C
$$
is an equivalence in homology with coefficients $\Z/q$. This
implies that $(P_{\Sigma B\Z/p} C)^{\wedge}_q \simeq
(BG)^{\wedge}_q$.

By R.~Levi's dichotomy result \cite[Theorem 1.1.4]{MR96c:55019}
for $q$-completions of classifying spaces of finite groups one
infers that $(BG)^{\wedge}_q$ must have infinitely many
non-trivial homotopy groups (because $G$ is $q$-perfect as its
abelianization must be $p$-torsion).
\end{proof}

\medskip

In view of the proof of the preceding theorem, the main question
is to determine whether or not the higher homotopy groups of the
cellularization of classifying spaces may contain $p$-torsion. We
make next an observation about $\Sigma B\Z/p$-local spaces which
play such an important role for the $B\Z/p$-cellularization
because of Chach\'olski's theorem~\ref{betterCW}. The proof is
very much in the spirit of the Dwyer-Wilkerson dichotomy result
\cite[Theorem~1.3]{MR92b:55004}, Levi's one
\cite[Theorem~1.1.4]{MR96c:55019}; compare also with Grodal's work
on Postnikov pieces \cite{MR1622342}.

\begin{proposition}
\label{sigmalocal}
Let $X$ be a simply-connected torsion $\Sigma B\Z/p$-local space.
Then it is either $p$-torsion free or has infinitely many homotopy
groups containing $p$-torsion.
\end{proposition}

\begin{proof}
Let $n$ be an integer $\geq 2$ and consider the Postnikov
fibration
$$
X \langle n \rangle \longrightarrow X \langle n-1 \rangle
\longrightarrow K(\pi_n X, n).
$$
Since the base point component of the iterated loop space
$\Omega^{n-1} X$ and $\Omega^{n-1} X \langle n-1 \rangle$ are
weakly equivalent there is a fibration
$$
\Omega^{n-1}_0 X
\longrightarrow K(\pi_n X, 1) \longrightarrow \Omega^{n-2} X
\langle n \rangle
$$
in which the fiber is $B\Z/p$-local. Thus the total space is
$B\Z/p$-local if and only if the base is so, i.e. the homotopy
group $\pi_n X$ contains $p$-torsion if and only if the
$n$-connected cover $X \langle n \rangle$ has some $p$-torsion. We
use here that a simply-connected $p'$-torsion space $Y$ has
trivial $p$-completion, therefore the pointed mapping space
$\map_*(B\Z/p, Y)$ is contractible.
\end{proof}

We turn now to a more detailed study of the case when $G$ is not a
$p$-group. We have seen that $CW_{B\Z/p} BG$ has infinitely many
non-trivial homotopy groups, but they can arise in two different
ways, because the space $P_{\Sigma B\Z/p} C$ may contain
$p$-torsion or not. If it is $p$-torsion free, it is weakly
equivalent to $P_{B\Z/p} BG \simeq \prod_{q \neq p}
BG^{\wedge}_q$. This means precisely that the cellularization
coincides with $\bar P_{B\Z/p} BG$, the homotopy fiber of the
nullification map.

\begin{theorem}
\label{step}
Let $G$ be a finite group. Then either the cellularization
$CW_{B\Z/p} BG$ has infinitely many homotopy groups containing
$p$-torsion or it fits in a fibration
$$
CW_{B\Z/p} BG \longrightarrow BG\longrightarrow \prod_{q\neq p}
(BG)^{\wedge}_q ,
$$
where the (finite) product is taken over all primes $q$ dividing
the order of $G$, and the right map is the product of the
completions.
\end{theorem}

\begin{proof}
As in the proof of Theorem~\ref{dichotomy} we may assume that $G$
is a $\Z/p$-cellular group. The homotopy cofiber of the evaluation
map $\vee_{[B\Z/p, BG]} B\Z/p \rightarrow BG$ is thus 1-connected
and so is $P_{\Sigma B\Z/p} C$. Moreover the group $G$ is finite,
so its integral homology groups are finite. The space $P_{\Sigma
B\Z/p} C$ is constructed out of $BG$ by taking iterated homotopy
cofibers of maps out of (finite) wedges of suspensions of $B\Z/p$,
whose homology groups are all torsion, and thus $P_{\Sigma B\Z/p}
C$ satisfies the conditions of the above proposition.

Assume thus that $P_{\Sigma B\Z/p} C$ is $p$-torsion free. As this
space is simply-connected, we see that it is actually
$B\Z/p$-local, and so $P_{B\Z/p} BG \simeq P_{\Sigma B\Z/p} C$.
The computation in \cite[Section 3]{Ramon} yields now that
$$
P_{B\Z/p} BG \simeq \prod_{q \neq p} BG^{\wedge}_q
$$
because $G$ is $\Z/p$-cellular (and thus equal to its
$\Z/p$-radical).
\end{proof}

\begin{example}
\label{S3}
{\rm Consider the $C_2$-cellular group $\Sigma_3$, symmetric group
on three letters. The choice of any transposition yields a map
$f:B\Z/2 \rightarrow B\Sigma_3$ which induces an isomorphism in
mod $2$ homology. Therefore the homotopy cofiber $C_f$ of $f$ is
$2$-torsion free and we are in case (i). This shows that
$CW_{B\Z/2} B\Sigma_3$ is a space whose fundamental group is
$\Sigma_3$ and its universal cover is $\Omega
(B\Sigma_3)^{\wedge}_3 \simeq S^3\{3\}$, the homotopy fiber of the
degree map on the sphere $S^3$, compare with \cite[VII,
4.1]{Bousfield72} and \cite[Theorem~7.5]{MR2003a:55024}.}
\end{example}

\section{The $p$-torsion free case}

The example of the symmetric group $\Sigma_3$ might lead us to
think that there will be very few cases when the universal cover
of $CW_{B\Z/p} BG$ is $p$-torsion free. In this section we will
see that there are surprisingly many groups for which this occurs
and one could even wonder if $p$-torsion actually can appear in
the $B\Z/p$-cellularization of $BG$.

Let us again consider the cofibration $\vee B\Z/p \rightarrow BG
\rightarrow C$ for a $\Z/p$-cellular group $G$. The cofiber $C$ is
a simply connected torsion space, hence equivalent to the finite
product $\prod C^{\wedge}_q$ where $q$ runs over the primes
dividing the order of $G$. We need to understand the $\Sigma
B\Z/p$-nullification of $C$. Since nullification commutes with
finite products and obviously $C^{\wedge}_q$ is $\Sigma
B\Z/p$-local for $q \neq p$, we only look at $C^{\wedge}_p$ and
infer that $CW_{B\Z/p} BG$ has infinitely many homotopy groups
with $p$-torsion if and only if $P_{\Sigma B\Z/p} C^{\wedge}_p$
does so, which happens if and only if $P_{\Sigma B\Z/p}
C^{\wedge}_p$ is not contractible by Theorem~\ref{step}.

\begin{proposition}
\label{nontrivialmap}
Let $G$ be a finite $\Z/p$-cellular group. For $CW_{B\Z/p} BG$ to
have infinitely many homotopy groups with $p$-torsion, there must
exist a $p$-complete space $Z$ which is $\Sigma B\Z/p$-local and a
map $f: BG \rightarrow Z$ which is not null-homotopic such that
the restriction $B\Z/p \rightarrow BG \rightarrow Z$ to any cyclic
subgroup of order $p$ is null-homotopic.
\end{proposition}

\begin{proof}
We need to understand when $P_{\Sigma B\Z/p} C^{\wedge}_p$ is not
contractible, or in other words when the map $C^{\wedge}_p
\rightarrow *$ is not a $\Sigma B\Z/p$-equivalence. This means by
definition that there exists some $\Sigma B\Z/p$-local space $Z$
for which the pointed mapping space $\map_*(C^{\wedge}_p, Z)$ is
not contractible. Because $C^{\wedge}_p$ is a  $p$-torsion space
(homotopy cofiber of such spaces) which is simply connected, we
can assume, using Sullivan's arithmetic square, that $Z$ is
$p$-complete.

The cofibration sequence $\vee B\Z/p \rightarrow BG \rightarrow C$
yields a fibration
$$
\map_*(C, Z) \rightarrow \map_*(BG, Z) \rightarrow \prod \map_*(B\Z/p, Z)
$$
in which the loop space of the base is trivial since $Z$ is
$\Sigma B\Z/p$-local. The existence of the map $f$ tells us that
there is a component of the total space (different from the
component of the constant map) which lies over the component of
the base point in the base. Therefore $\map_*(C, Z)$ has at least
two component as well. As $Z$ is $p$-complete $\map_*(C, Z)$ is
weakly equivalent to $\map_*(C^{\wedge}_p, Z)$.

On the other hand, if no such map $f$ exists, $\map_*(C, Z)$ is
weakly equivalent to the connected component of the constant
$\map_*(BG, Z)_c$, which is contractible by Dwyer's result
\cite[Theorem~1.2]{MR97i:55028}.
\end{proof}

Our aim is now to explain why the map $f$ in the above proposition
cannot exist in many cases.

\begin{proposition}
\label{noptorsion}
Let $G$ be a finite $\Z/p$-cellular group whose Sylow $p$-subgroup
is generated by elements of order $p$. The universal cover of
$CW_{B\Z/p} BG$ is then $p$-torsion free.
\end{proposition}

\begin{proof}
By Theorem~\ref{step} the universal cover of $CW_{B\Z/p} BG$ is
either $p$-torsion free or it has infinitely many homotopy groups
with $p$-torsion. As $CW_{B\Z/p} BG$ is the homotopy fiber of a
map $BG \rightarrow P_{\Sigma B\Z/p} C$, the same is true for
$P_{\Sigma B\Z/p} C^{\wedge}_p$. We prove now that this space is
always contractible because the assumptions on the map $f$ in
Proposition~\ref{nontrivialmap} are never satisfied.

Let $f: BG \rightarrow Z$ be a map into a $\Sigma B\Z/p$-local
space which is trivial when restricted to any cyclic subgroup
$\Z/p$ in $G$. The composite $BG \rightarrow Z \rightarrow K(\pi_1
Z, 1)$ is then null-homotopic because it corresponds to a
homomorphism $G \rightarrow \pi_1 Z$ restricting trivially to all
generators of $G$. Therefore the map $f$ lifts to the universal
cover of $Z$ and so we might assume that $Z$ is $1$-connected. We
can also $p$-complete it if necessary, so that $Z$ is actually
$H\Z/p$-local. By Dwyer's theorem \cite[Theorem~1.4]{MR97i:55028},
the map $f$ is null-homotopic if and only if the restriction $\bar
f$ to some Sylow $p$-subgroup $S$ of $G$ is so. Because we assume
$S$ is generated by elements of order $p$, we can consider as in
\cite[Proposition 4.14]{Ramon} the split extension $T
\hookrightarrow S \epi \Z/p$ given by a maximal normal subgroup
generated by all generators but one. By induction on the order the
composite
$$
BT \rightarrow BS \stackrel{\bar f}{\longrightarrow} Z
$$
is null-homotopic. Hence using Zabrodsky's Lemma (see
\cite[Proposition~3.5]{MR97i:55028} and \cite{CCS2} for a more
detailed account in our setting) we see that the map $\bar f$
factors through a map $g: B\Z/p \rightarrow Z$ and it is
null-homotopic if and only if $g$ is so. As we suppose that the
map $f$ is trivial when restricted to any cyclic subgroup $\Z/p$
we can conclude since the above extension is split.
\end{proof}

\begin{example}
\label{symmetric}
{\rm Consider the $C_2$-cellular group $\Sigma_{2^n}$, symmetric
group on $2^n$ letters with $n \geq 2$. The Sylow $2$-subgroup is
an iterated wreath product of copies of $\Z/2$ which is always
generated by elements of order $2$. Therefore the above
proposition applies and we obtain that $P_{\Sigma B\Z/2} C \simeq
\prod_{q \neq 2} (B\Sigma_{2^n} )^{\wedge}_q$. Thus $CW_{B\Z/2}
B\Sigma_{2^n}$ fits into a fibration
$$
 \prod_{q \neq 2} \Omega (B\Sigma_{2^n} )^{\wedge}_q \rightarrow CW_{B\Z/2} B\Sigma_{2^n}
 \rightarrow B\Sigma_{2^n}.
$$}
\end{example}

\medskip

In the next proposition we will see that the existence of
$p$-torsion in the upper homotopy of $CW_{B\Z/p} BG$ is strongly
related with the $B\Z/p$-cellularity of $BG^{\wedge}_p$. Recall
that the fundamental group of $BG^{\wedge}_p$ is always isomorphic
to the group theoretical $p$-completion $G^{\wedge}_p$, i.e. the
quotient of $G$ by $O^p(G)$, the maximal $p$-perfect subgroup
of~$G$.

\begin{proposition}
\label{completion}
Let $G$ be a finite group generated by order $p$ elements which is
not a $p$-group. Then the universal cover of $CW_{B\Z/p} BG$ is
$p$-torsion free if and only if the $p$-completion of $BG$ is
$B\Z/p$-cellular.
\end{proposition}

\begin{proof}
We know from \cite[Proposition 2.1]{Broto03} that the completion
map $BG\longrightarrow BG^{\wedge}_p$ induces a chain of
bijections identifying the set of unpointed homotopy classes
$$
[B\Z/p,BG] \simeq [B\Z/p,BG^{\wedge}_p] \simeq
\textrm{Rep}(\Z/p,G).
$$
Choose a set of representatives $f: B\Z/p \rightarrow BG$ for all
conjugacy classes of elements of order $p$ in $G$ and write
$f^{\wedge}_p$ for the corresponding map into the $p$-completion
of $BG$. We have seen in Lemma~\ref{betterCW} that $CW_{B\Z/p} BG$
can be constructed as the homotopy fiber of the composite $BG
\rightarrow D \rightarrow P_{\Sigma B\Z/p} D$, where $D$ is the
homotopy cofiber of the evaluation map $ev: \bigvee B\Z/p
\rightarrow BG$. Likewise $CW_{B\Z/p} (BG)^{\wedge}_p$ is
constructed using the completed version, call the corresponding
cofiber $D'$. In short we have a diagram of cofibrations
$$
\xymatrix{ \bigvee B\Z/p \ar[r]^{\vee f} \ar@{=}[d] & BG \ar[r]
\ar[d] &
D \ar[d] \\
\bigvee B\Z/p \ar[r]^{\vee f^{\wedge}_p} & BG^{\wedge}_p \ar[r] &
D' }
$$
Now since $G$ (and thus its quotient $G^{\wedge}_p$ as well) is
generated by elements of order $p$, the cofibers $D$ and $D'$ are
simply connected. As $\bigvee B\Z/p$ and $BG$ are rationally
trivial, we use Sullivan's arithmetic square to establish that
$D=\prod_{q\textrm{ prime}}D^{\wedge}_q$ and since moreover
$BG^{\wedge}_p$ is $p$-complete, $D'$ is so.

Comparing the Mayer-Vietoris sequences in mod $p$ homology of the
two cofibrations, it is easy to see that $D^{\wedge}_p\simeq D'$.
Therefore we can compute
$$
P_{\Sigma B\Z/p} D \simeq \prod_q P_{\Sigma B\Z/p} D^{\wedge}_q
\simeq \prod_{q \neq p} D^{\wedge}_q \times P_{\Sigma B\Z/p} D' .
$$
Hence the universal cover of the $B\Z/p$-cellularization of $BG$
is $p$-torsion free if and only if $D'$ is $\Sigma B\Z/p$-acyclic.
This is equivalent for $BG^{\wedge}_p$ to be $B\Z/p$-cellular.
\end{proof}

Let us finish this section by showing with some examples the
applicability of the last result.

\begin{example}
\label{symmetric2}
{\rm Consider the $2$-completion of the classifying space of the
symmetric group $\Sigma_{2^n}$. According to \ref{symmetric}, the
Sylow 2-subgroup of $\Sigma_{2^n}$ is generated by order $2$
elements, and moreover $\pi_k(CW_{B\Z /2}B\Sigma_{2^n})$ is
2-torsion free if $k\geq 2$. Thus by Proposition~\ref{completion}
$(B\Sigma_{2^n})^{\wedge}_2$ is $B\Z /2$-cellular.}
\end{example}

Unlike to what happens with the cellularization of $BG$, our study
does not give in general information about
$CW_{B\Z/p}(BG^{\wedge}_p)$ when $G$ is not generated by order $p$
elements. In this case, however, it is sometimes possible to
reduce the problem to the case of groups for which the hypothesis
of the theorem hold. To show this, we compute one last example.

\begin{example}
\label{alternating}
{\rm Let us compute the $B\Z /2$-cellularization of
$(BA_4)^{\wedge}_2$. It is known that the natural inclusion
$A_4<A_5$ induces an equivalence in mod 2 homology, and then
$(BA_4)^{\wedge}_2\simeq (BA_5)^{\wedge}_2$. Now, as $A_5$ is
simple with 2-torsion, it is generated by order $2$ elements, and
moreover the Sylow 2-subgroup of $A_5$ is the Klein group $\Z
/2\times\Z /2$. Hence, $(BA_5)^{\wedge}_2$ is $B\Z /2$-cellular by
Proposition~\ref{completion}, and then $(BA_4)^{\wedge}_2$ is so
as well.}
\end{example}

\section{Cellularization and fusion}

To find an example where $CW_{B\Z/p} BG$ has $p$-torsion, one
should look by Proposition~\ref{noptorsion} for groups generated
by elements of order $p$ with a Sylow $p$-subgroup which is not
generated by elements of order $p$. We will see in this section
that this is by far not a sufficient condition, as illustrated
below by the example of the the finite simple group $PSL_3(3)$. We
start with a simple observation which leads then naturally to a
closer analysis of the fusion of the groups we look at.

\begin{proposition}
\label{fusionsubgroup}
Let $G$ be a group generated by elements of order $p$ such that
its Sylow $p$-subgroup $S$ is not so. If $S$ is generated by the
elements of $\Omega_1(S)$ together with all their conjugates by
elements of $G$ which belong to $S$ then $CW_{B\Z/p} BG$ is
$p$-torsion free.
\end{proposition}

\begin{proof}
We have seen in Proposition~\ref{nontrivialmap} that the existence
of $p$-torsion is detected by a non-trivial map $f: BG^{\wedge}_p
\rightarrow Z$ into some $\Sigma B\Z/p$-local space $Z$. Such a
map is null-homotopic if and only if the restriction to $BS$ is
so. Pick a generator in $S$ and consider the cyclic subgroup $T$
it generates. By assumption this subgroup is conjugate in $G$ to
some subgroup of $\Omega_1(S)$. Because conjugation in $G$ induces
a weak equivalence on $BG$ we see that $f$ is null-homotopic when
restricted to $BT$ (it is so when restricted to $B\Omega_1(S)$).
An induction on the order of the Sylow subgroup as in the proof of
Proposition~\ref{noptorsion} allows then to conclude that $f$
itself cannot be essential.
\end{proof}

\begin{example}
\label{Thevenaz}
{\rm The symmetric group $\Sigma_3$ acts by permutation on
$(\Z/4)^3$. The diagonal is invariant, and so is the ``orthogonal"
subgroup $\Z/4 \times \Z/4$. We define $G$ to be the semi-direct
product of $\Z/4 \times \Z/4$ by $\Sigma_3$. It is easy to check
that $G$ is generated by elements of order $2$, but the Sylow
2-subgroup $S = (\Z/4 \times \Z/4) \rtimes \Z/2$ is not. The
subgroup $\Omega_1(S)$ has index 2, and a representative of the
generator of the quotient can be taken inside $S$ to have order
four (it is inside of $G$ a product of three elements of order 2).
This element has a conjugate which lies inside $\Omega_1(S)$, so
we may conclude by the above proposition that $CW_{B\Z/2} BG$ is
$2$-torsion free, i.e. its universal cover is $\Omega
(BG^{\wedge}_3)$ by Theorem~\ref{step}. }
\end{example}

\begin{example}
\label{L33}
{\rm The group $G= PSL_3(3)$ presents the same features as the
above example. It is generated by elements of order 2 (it is
simple), but its Sylow 2-subgroup $S$ is semi-dihedral of order
16. The subgroup $\Omega_1(S)$ is dihedral of order $8$. Here as
well the universal cover of $CW_{B\Z/2} BPSL_3(3)$ is $\Omega
\bigl( BPSL_3(3)^{\wedge}_3\bigr)$.}
\end{example}

The previous examples share the common property that the index of
the subgroup $\Omega_1(S)$ inside the Sylow $2$-subgroup is $2$.
In this case a non-trivial map $BG^{\wedge}_2 \rightarrow Z$ as
desired can never exist. This comes from the fact that the only
automorphism of the quotient $\Z/2 \cong \Omega_1(S)/S$ is the
identity, and thus such a map should actually factor through
$B\Z/2$, which is obviously impossible. Hence, we should look for
groups where the subgroup $\Omega_1(S)$ has large index in $S$,
and is preserved by fusion. The Suzuki group $Sz(2^n)$ with $n$ an
odd integer $\geq 3$ is such a group, as we learned from Bob
Oliver. The section \cite[16.4]{Gor80} is extremely useful to
understand the basic subgroup and fusion properties of the Suzuki
groups. In particular the normalizer of the Sylow $2$-subgroup of
$Sz(2^n)$ is a semi-direct product $S \rtimes \Z/(2^n-1)$ which is
maximal in $Sz(2^n)$.

\begin{lemma}
\label{upto}
Let $S$ denote the Sylow $2$-subgroup of $Sz(2^n)$. The inclusion
of the maximal subgroup $S \rtimes \Z/(2^n-1) \hookrightarrow
Sz(2^n)$ induces a weak equivalence of $2$-complete spaces
$B\bigl(S \rtimes \Z/(2^n-1)\bigr)^{\wedge}_2 \simeq
BSz(2^n)^{\wedge}_2$.
\end{lemma}

\begin{proof}
The 2-Sylow subgroup of $Sz(2^n)$ can be written as an extension
$(\Z /2)^n\hookrightarrow S\twoheadrightarrow (\Z /2)^n$, where
the kernel is the center of the group and contains all its order 2
elements. Observe in particular that $S$ is \emph{not} generated
by order 2 elements, an unavoidable condition according to
\ref{noptorsion}, and its index is $2^n$.

In order to understand the $2$-completion of $BSz(2^n)$ we have
learned from \cite{Broto032} that we need to determine which
subgroup of $S$ are 2-centric and 2-radical (see also
\cite[2.1]{Jackson04} for a definition). Our group $Sz(2^n)$ has
the property that all of its 2-Sylow subgroups are disjoint (they
have trivial intersection), hence the unique 2-centric 2-radical
subgroup (up to conjugation) is $S$ itself. The outer automorphism
group of $S$ in the fusion system of the Suzuki group, i.e. the
quotient of the normalizer $N_{Sz(2^n)}(S)$ by $S$, is cyclic of
order $2^n-1$, generated by an element $\phi$ which acts
fixed-point free and permutes transitively the non-trivial
elements of the center of $S$ (see \cite[16.4]{Gor80}).

Now it is clear from the construction that the inclusion
$S\rtimes\Z /(2^n-1) \hookrightarrow Sz(2^n)$ induces an
isomorphism of fusion and linking systems in the sense of
\cite{Broto032}, and in particular it induces a homotopy
equivalence $f:B(S\rtimes\Z /(2^n-1))^{\wedge}_2\simeq
BSz(2^n)^{\wedge}_2$.
\end{proof}

\begin{lemma}
\label{representation}
Consider the semi-direct product $(\Z/2)^n \rtimes \Z/(2^n-1)$
where a generator $\phi$ of $\Z/(2^n-1)$ acts on the elementary
abelian group of rank $n$ by permuting transitively the $2^n-1$
non-trivial elements. There exists then a faithful representation
$\sigma: (\Z/2)^n \rtimes \Z/(2^n-1) \hookrightarrow U(2^n-1)$.
\end{lemma}

\begin{proof}
The representation can be induced from the trivial one on the
subgroup $(\Z/2)^n$. When $n=3$ we can be very explicit. Send the
first generator of the elementary abelian group to the diagonal
matrix with entries $(-1, 1, -1, -1, -1, 1, 1)$, and the other
elements to the cyclic permutations of it. The standard cyclic
permutation matrix of order $7$ in $U(7)$ is the image of $\phi$
and so $\sigma$ is well-defined.
\end{proof}

\begin{proposition}
\label{Suzuki}
There exists a non-trivial map $BSz(2^n)\longrightarrow
BU(2^n-1)^{\wedge}_2$ such that the composition $B\Z /2
\longrightarrow BSz(2^n)\longrightarrow BU(2^n-1)^{\wedge}_2$ is
null-homotopic for every cyclic subgroup $\Z /2$ in $Sz(2^n)$.
\end{proposition}

\begin{proof}
We construct actually a map from the $2$-completion of $BSz(2^n)$
and the desired map is then obtained by pre-composing with the
completion map $BSz(2^n) \rightarrow BSz(2^n)^{\wedge}_2$. By
Lemma~\ref{upto} we only need to construct a map out of $B \bigl(
(S \rtimes \Z/(2^n-1) \bigr)^{\wedge}_2$ where $S$ denotes the
2-Sylow subgroup of $Sz(2^n)$.

Because $BU(2^n-1)$ is simply connected the sets of homotopy
classes of pointed and unpointed maps into $BU(2^n-1)$ agree. The
fusion system of the Suzuki group (i.e. that of the semi-direct
product) is reduced to the Sylow subgroup, so that the set of
$[BSz(2^n), BU(2^n-1)]$ is isomorphic to the set of fusion
preserving representations of $S$ inside $U(2^n-1)$. This means
that for finding a non-trivial map $f: BSz(2^n) \rightarrow
BU(2^n-1)^{\wedge}_2$ it is enough to find a representation $\rho:
S \rtimes \Z/(2^n-1) \rightarrow U(2^n-1)$ which is non-trivial
when restricted to $S$ (compare with Dwyer's result
\cite[Theorem~1.4]{MR97i:55028} which we have already used
before).

On the other hand we want the map $f$ to be null-homotopic when
restricted to any cyclic subgroup $\Z/2$. In other words the
composite $\Omega_1(S) \rightarrow S \rtimes \Z/(2^n-1)
\rightarrow U(2^n-1)$ should be trivial. The subgroup
$\Omega_1(S)$ generated by all elements of order $2$ is its
center, an elementary abelian subgroup of rank $n$. It is normal
in $S \rtimes\Z /(2^n-1)$, and the quotient is isomorphic to
$(\Z/2)^n \rtimes\Z /(2^n-1)$, where the action of the generator
$\phi$ of order $2^n-1$ permutes transitively the non-trivial
elements (one has a bijection between these non-trivial classes in
the quotient and the squares of their representatives in the
center of $S$). Hence we can define $\rho$ as the composite
$$
S \rtimes\Z /(2^n-1) \longrightarrow (\Z/2)^n \rtimes\Z /(2^n-1)
\stackrel{\sigma}{\longrightarrow} U(2^n-1)
$$
where $\sigma$ is the representation constructed in the preceding
lemma. It is convenient to remark here that the existence of such
a representation does not contradict the theorems in this article
since $S \rtimes \Z/(2^n-1)$ is \emph{not} generated by order 2
elements.

Observe that the morphism $B(S \rtimes \Z/(2^n-1))^{\wedge}_2
\longrightarrow BU(2^n-1)^{\wedge}_2$ induced by the
representation we have just constructed is clearly trivial when
pre-composing with any map $B\Z /2 \longrightarrow B\bigl(
S\rtimes\Z /(2^n-1) \bigr)^{\wedge}_2$, because the (unpointed)
homotopy classes of these last maps can be identified with the
conjugacy classes of $\mathbb{Z}/2$ inside $S\rtimes\Z /(2^n-1)$.
\end{proof}

\begin{remark}
{\rm Note that the representation $\rho$ constructed in the
previous proposition cannot be induced by a homomorphism
$Sz(2^n)\longrightarrow U(2^n-1)$, because the group $Sz(2^n)$ is
simple, hence generated by order 2 elements. If it were so, the
homomorphism would be zero over the generators, and thus trivial.}
\end{remark}

The methods above can also be used to obtain examples at odd
primes, as we sketch in the sequel. The following example was
pointed out to us by Antonio Viruel.

\begin{example}
\label{oddprimes}
{\rm Let $p$ be an odd prime, $n$ an integer $\geq 2$, and
$q=mp^n+1$. Consider the linear group $PSL_2(q)$. According to
(\cite[15.1.1]{Gor80}), the $p$-Sylow subgroup of $PSL_2(q)$ is
cyclic of order $p^n$, and moreover its normalizer
$N_{PSL_2(q)}(S)$ is isomorphic to the semidirect product
$\Z/{p^n}\rtimes\Z /2$, where the action is given by the change of
sign. As the Sylow subgroup is abelian, it can be deduced from
(\cite[3.4]{DRV04}) that the inclusion
$N_{PSL_2(q)}(S)\hookrightarrow PSL_2(q)$ induces a homotopy
equivalence between the 2-completions of the classifying spaces.
Reasoning as above, it is enough to give, for some $k$, a
non-trivial map $\Z/{p^n}\rtimes\Z /2\longrightarrow U(k)$ which
is trivial when restricted to order $p$ elements. We first
construct an inclusion $j: \Z/p\rtimes\Z /2 \hookrightarrow U(2)$
as follows. If $x = e^{2\pi i /p}$  the image of the generator of
order $p$ is the diagonal matrix with entries $(x, \bar x)$, and
the image of the generator of order $2$ is the standard
permutation matrix. One can then take for example the composition
$$
f:\Z/{p^n}\rtimes\Z /2\longrightarrow\Z/p\rtimes\Z
/2\stackrel{j}{\longrightarrow} U(2),
$$
where the first map is the natural projection. Now the induced map
at the level of $p$-completed classifying spaces
$$
(BPSL_2(q))^{\wedge}_p\simeq (B(\Z/{p^n}\rtimes\Z /2))^{\wedge}_p\longrightarrow
BU(2)^{\wedge}_p
$$
is essential, but homotopically trivial when
precomposing with any map from $B\Z /p$.}
\end{example}

We are finally able to show that the $B\Z/p$-cellularization of
$BG$'s does not always coincide with the fiber of the
$B\Z/p$-nullification. Our example is of course given by the
groups described above. This implies that our cellular dichotomy
Theorem~\ref{dichotomy} is in fact a trichotomy result: the higher
homotopy groups of $CW_{B\Z/p} BG$ are either all trivial, or
$p$-torsion free and infinitely many are non-trivial, or else
infinitely many do contain $p$-torsion.

\begin{theorem}
\label{SuzukiCW}
For every integer $n$ of the form $2^{2k+1}$ the
$B\Z /2$-cellularization of $BSz(n)$ has 2-torsion in an infinite
number of homotopy groups. Likewise, if $p$ is an odd prime and
$q$ is any integer of the form $mp^k+1$ with $k \geq 2$, then
the $B\Z /p$-cellularization of  $BPSL(q)$ has $p$- torsion in an infinite
number of homotopy groups.
\end{theorem}

\begin{proof}
The loop space $\Omega BU(2^n-1)^{\wedge}_2$ is $B\Z/2$-local by
\cite[9.9]{Miller}. It is now a direct consequence of
Proposition~\ref{nontrivialmap} that the map $f: BSz(2^n)
\rightarrow BU(2^n-1)^{\wedge}_2$ constructed in the previous
proposition implies the existence of infinitely many homotopy
groups $\pi_n BSz(2^n)$ containing $2$-torsion. The same argument
can be used at odd primes for the linear groups described in
example \ref{oddprimes}.
\end{proof}


\bibliographystyle{alpha}
\bibliography{dichotomy}

\bigskip

\bigskip\noindent
Ram\'on J. Flores, J\'er\^ome Scherer

\noindent Departament de Matem\`atiques, Universitat Aut\'onoma de
Barcelona, E--08193 Bellaterra \\
e-mail: {\tt ramonj@mat.uab.es, jscherer@mat.uab.es}

\end{document}